\documentclass[11pt,twoside]{article}
 \usepackage{a4}
\usepackage{amssymb}
\usepackage{amsmath}

\newtheorem{theorem}{Theorem}
\newtheorem{FubiniTheorem}{Fubini's Theorem}

\newtheorem{conjecture}{Conjecture}

\newtheorem{lemma}{Lemma}

\renewcommand{\Bbb}[1]{\mathbb{#1}}
\newcommand{\N}{{\Bbb N}}         
\newcommand{\R}{{\Bbb R}}         
\newcommand{\Rp}{{\Bbb R}_+}    
\newcommand{\Z}{{\Bbb Z}}         

\newcommand{\cD}{{\cal D}}

\newcommand{\cS}{{\cal S}}

\newcommand{\ve}{\varepsilon}

\newcommand{\vv}[1]{{\mathbf{#1}}}

\renewcommand{\le}{\leqslant}
\renewcommand{\ge}{\geqslant}
\newcommand{\nz}{\smallsetminus\{\vv 0\}}

\begin{document}

\title{Multiplicative zero-one laws and  \\ metric number theory}

\author{Victor Beresnevich\footnote{EPSRC Advanced Research Fellow, grant EP/C54076X/1}
\and
Alan Haynes\footnote{Research supported by EPSRC grant EP/F027028/1}
\and Sanju Velani\footnote{Research supported by EPSRC grants EP/E061613/1 and EP/F027028/1 }
}

\date{}

\maketitle

\vspace*{-4ex}

\begin{abstract}
We develop the classical theory of Diophantine approximation without assuming monotonicity or convexity. A complete `multiplicative' zero-one law is established akin to the `simultaneous' zero-one laws of Cassels and Gallagher.  As a consequence we are able to establish the analogue of the Duffin-Schaeffer theorem within the multiplicative setup.  The key ingredient is the rather simple but nevertheless  versatile `cross fibering principle'. In a nutshell it enables us to `lift' zero-one laws to higher dimensions.
\end{abstract}

\textit{Keywords}: Zero-one laws, metric multiplicative Diophantine approximation, Duffin-Schaeffer theorem

\textit{Subject classification}: 11J13, 11J83, 11K60

\section{Introduction}

The  theory of multiplicative Diophantine approximation
is concerned with the set
$$
\cS^\times_n(\psi):=
\{(x_1,\dots,x_n)\in [0,1]^n:\prod_{i=1}^n\|qx_i\|<\psi(q)\text{ for i.m. }q\in\N\},
$$
where $\|qx\|=\min\{|qx-p|:p\in\Z\}$, `i.m.' means `infinitely many' and $\psi:\N\to\Rp$  is a a non-negative  function.  For obvious reasons the function $\psi$  is often referred to as an approximating function. For convenience, we work within the unit cube $[0,1]^n$ rather than $\R^n$; it makes full measure results easier to state and avoids ambiguity.  In fact, this is not at all restrictive since  the set under consideration  is invariant under translation by integer vectors.

Multiplicative Diophantine approximation is currently an active area of research. In particular, the long standing conjecture of Littlewood that states that $\cS^\times_2(q\mapsto \ve q^{-1})=\R$ for any $\ve>0$ has attracted much attention -- see  \cite{badvel,pvl,vent} and references within. In this paper we will address the multiplicative analogue of yet another long standing classical problem; namely, the Duffin-Schaeffer conjecture.

Given $q\in\N$ and $x\in\R$, let $$\|qx\|':=\min\{|qx-p|:p\in\Z,\ (p,q)=1\}  \, , $$ and consider  the standard simultaneous  sets
$$
\cD_n(\psi):=
\{(x_1,\dots,x_n)\in [0,1]^n:
\big(\max\limits_{1\le i \le n}\|qx_i\|'\big)^n<\psi(q)\text{ for i.m. }q\in\N\}  \,
$$
and
$$
\cS_n(\psi):=
\{(x_1,\dots,x_n)\in [0,1]^n:
\big(\max\limits_{1\le i \le n}\|qx_i\|\big)^n<\psi(q)\text{ for i.m. }q\in\N\}  \, .
$$

\noindent An elegant measure theoretic property of these sets is that they are always of zero or full Lebesgue measure $ | \ . \ | $ irrespective of the dimension or the approximating function.   Formally, for $n\ge 1$ and any non-negative function  $\psi:\N\to\Rp$
\begin{equation}\label{zo}
|\cS_n(\psi)|\in\{0,1\}  \qquad  {\rm and}    \qquad|\cD_n (\psi)|\in\{0,1\}  \, .
\end{equation}

\noindent
The former  zero-one law is due to Cassels \cite{Cassels-50:MR0036787} while the latter is due to  Gallagher \cite{Gallagher-61:MR0133297} when  $n=1$ and Vilchinski \cite{vil} for  $n$ arbitrary.   By making use of a refined version  of Cassels' zero-one law,
Gallagher \cite{Gallagher65} proved that for $n\ge 2$
\begin{equation}\label{gal}
  |\cS_n(\psi)|= 1  \qquad \text{if } \quad   \sum_{q=1}^\infty \psi(q) = \infty \; .
\end{equation}

\noindent{\em Remark 1.\ } Regarding the above statement and indeed the statements and  conjectures below,
by making use of the  Borel-Cantelli Lemma from probability theory, it is straightforward to establish the complementary convergent results;
i.e. if the sum in question converges then the set in question is of zero measure.

\medskip

\noindent The case that $n=1$ is excluded from the statement given by (\ref{gal}) since it is false. Indeed, Duffin \& Schaeffer \cite{DuffinSchaeffer} gave a counterexample and formulated an alternative appropriate statement.
The Duffin-Schaeffer conjecture\footnote{To be precise  Duffin and Schaeffer stated their conjecture for $n=1$. The higher dimensional version is attributed  to Sprind\v{z}uk -- see \cite[pg63]{Sprindzuk-1979-Metrical-theory}.}  states that
\begin{equation}\label{ds1}
  |\cD_n(\psi)|= 1  \qquad \text{if } \quad   \sum_{q=1}^\infty\left(\frac{\varphi(q)}{q}\right)^n\psi(q) = \infty \; ,
\end{equation}
where $\varphi$ is the Euler phi function. The consequence of the  zero-one law for $ \cD_n(\psi)$ is that it reduces the Duffin-Schaeffer conjecture to  showing that $|\cD_n(\psi)|>0$. Using this fact the conjecture has been established in the case $n\ge 2$ by Pollington $\&$ Vaughan \cite{Pollington-Vaughan-1990}.
Although various partial results have been obtained in the case $n=1$, the full
conjecture represents a key unsolved problem in number
theory.  For background and recent developments regarding this fundamental problem see \cite{Beresnevich-Bernik-Dodson-Velani-Roth,DuffinSchaeffer,Harman-1998a,haypollvel}.
However, it is worth  highlighting the Duffin-Schaeffer theorem which states that (\ref{ds1}) holds whenever
 $$
 \limsup_{Q\to\infty}\left(\sum_{q=1}^{ Q} \left(\frac{\varphi (q)}{q}\right)\psi (q)\right)  \left(\sum_{q=1}^{ Q}  \psi (q) \right)^{-1} \, >  \, 0  \ .
$$
Note that this condition implies that the convergence/divergence properties of the sums in (\ref{gal}) and (\ref{ds1}) are equivalent.

As already mentioned, the purpose of this paper is to  consider the multiplicative setup and  in particular, the multiplicative analogue of the Duffin-Schaeffer conjecture. With this in mind, it is natural to define the set
$$
\cD^\times_n(\psi):=
\{(x_1,\dots,x_n)\in [0,1]^n:
\prod_{i=1}^n\|qx_i\|'<\psi(q)\text{ for i.m. }q\in\N\}.
$$

\

The ultimate goal  is to prove the following two statements.

\begin{conjecture} \label{con1}
Let $n\ge 2$ and $\psi:\N\to[0,\frac12)$. Then
\begin{equation}\label{c1}
|\cS^\times_n(\psi)|=  1 \qquad  \text{if } \quad \sum_{q=1}^\infty \psi(q)\left(\log \psi(q)^{-1}\right)^{n-1}=\infty  \, ,
\end{equation}
\end{conjecture}


\begin{conjecture} \label{con2}
Let $n\ge 1$ and $\psi:\N\to[0,\frac12)$. Then
\begin{equation}\label{c2}
|\cD^\times_n(\psi)|= 1  \qquad \text{if } \quad \sum_{q=1}^\infty\left(\frac{\varphi(q)}{q}\right)^n \psi(q)\left(\log \psi(q)^{-1}\right)^{n-1}=\infty \, .
\end{equation}
\end{conjecture}

\noindent   Throughout the paper,
$$\psi(q)\left(\log \psi(q)^{-1}\right)^{n-1} :=0   \quad {\rm  whenever} \quad  \psi(q)=0 \, . $$
 In view of the Duffin-Schaeffer counterexample it is necessary to exclude $n=1$ from the statement of Conjecture~\ref{con1}. Clearly,  the Duffin-Schaeffer conjecture and Conjecture~\ref{con2} coincide when $n=1$.

\medskip

\noindent{\em Remark 2.\ } For $n \ge 2$, the results of Gallagher and  Pollington $\&$ Vaughan  establish the analogues of the above conjectures for the standard simultaneous sets $\cS_n(\psi)$ and $\cD_n(\psi)$.

\subsection{The story so far: convexity versus monotonicity  \label{conmom}}
Throughout this section, assume that $n \ge 2$.
Geometrically, the multiplicative sets $\cS_n^\times(\psi)$ and $\cD_n^\times(\psi)$  consist of points in the  unit cube that lie within infinitely many `hyperbolic' domains
$$
{\rm H} = {\rm H}(\psi,\vv p,q) :=  \Big\{ \vv x \in [0,1]^n:
\textstyle{\prod_{i=1}^n |x_i -p_i/q |<\psi(q)/q^n}
\Big\}
$$
centered around rational points $\vv p /q$ where $\vv p =(p_1, \ldots, p_n) \in \Z^n$ and $ q \in  \N$.  In the case of $ \cD_n^\times(\psi)$ we impose the additional co-primeness condition $(p_i,q)=1$  on the rational points.  The approximating function $\psi$ governs the  size of the domains ${\rm H}$.
In the case of the standard simultaneous sets $\cS_n(\psi)$ and $\cD_n(\psi)$ the domains ${\rm H}$ are replaced by the `cubical' domains
$$ {\rm C} = {\rm C}(\psi,\vv p,q) :=  \Big\{ \vv x\in [0,1]^n:
\textstyle{\max\limits_{1\le i \le n}|x_i -p_i/q|^n<\psi(q)/q^n}
\Big\}.
$$
The significant difference between the standard and multiplicative situation  is that the domains ${\rm C}$  are convex while the domains ${\rm H}$ are non-convex.  It is this difference that lies behind the fact that Conjectures~\ref{con1} \& \ref{con2} are still open whilst their standard simultaneous counterparts have been established --  recall we assume that $n \geq 2$. In short, without imposing additional assumptions,  convexity is vital in the methods employed by Gallagher and  Pollington $\&$ Vaughan  to establish  \eqref{gal} and  \eqref{ds1} respectively. Indeed, their methods can be refined and adapted to deal with $\limsup$ sets arising from more general convex domains  but convexity itself  seems to be unremovable -- see \cite[Chp.3]{Harman-1998a} and references within.  However, the landscape is completely different if we impose the additional assumption that the approximating function $\psi$ is monotonic.  For instance we can then  overcome the fact that the domains ${\rm H}$  associated with the sets $\cS_n^\times(\psi)$ and $\cD_n^\times(\psi)$  are  non-convex and Conjectures~\ref{con1}~\&~\ref{con2} correspond to a well known theorem of Gallagher \cite{gallagher1962}. In fact,  Gallagher considers  $\limsup$ sets arising from more general  domains but monotonicity plays a crucial role in his approach and seems to be unremovable.   Note that for monotonic $\psi$ the convergence/divergence properties of the sums appearing in \eqref{c1} and \eqref{c2}  are equivalent and since   $\cS_n^\times(\psi)  \supset \cD_n^\times(\psi)$  it follows that  Conjecture~\ref{con2} implies Conjecture~\ref{con1}.

The upshot is that the current body of metrical results for $\limsup$ sets requires  that either the approximating domains are  convex  or that the approximating function is  monotonic.


\subsection{Statement of results \label{state}}
Our first theorem is the multiplicative analogue of the Cassels-Gallagher zero-one law. It reduces Conjectures~\ref{con1} \& \ref{con2} to showing that the corresponding sets are of positive measure. In principal,  it is easier to prove positive measure statements than full measure statements. More to the point,  there is a well established  mechanism in place  to obtain  lower bounds for the measure of $\limsup$ sets -- see \S\ref{ohya} below or \cite[\S8]{Beresnevich-Dickinson-Velani-06:MR2184760} for a more comprehensive account.

\begin{theorem}\label{t1}  Let $n\ge 1$ and $\psi:\N\to\Rp$. Then
$$|\cS_n^\times(\psi)|\in\{0,1\}  \qquad  and   \qquad |\cD_n^\times(\psi)|\in\{0,1\}  \, . $$
\end{theorem}

\noindent The proof will rely on the general  technique developed in \S\ref{cross} which we refer to as the  {\em cross fibering principle}.
Given its simplicity, we suspect that it may well have applications elsewhere in one form  or another.

\medskip

The following theorem represents our `direct' contribution to Conjectures \ref{con1} \& \ref{con2} and is the complete multiplicative analogue of the Duffin-Schaeffer theorem.

\begin{theorem}\label{MDST}
Let $n\ge 1$, $\psi:\N\to[0,\frac12)$.
Then
$$|\cS_n^\times(\psi)|  = 1 =  | \cD_n^\times(\psi)|
\qquad \text{if } \qquad  \sum_{q=1}^{\infty} \psi(q)\left(\log \psi(q)^{-1}\right)^{n-1}
  = \infty
$$
and
\begin{equation}\label{cond1}
 \limsup_{Q\to\infty}\,{\displaystyle\sum_{q=1}^{ Q} \left(\frac{\varphi(q)}{q}\right)^n \psi(q)\left(\log \psi(q)^{-1}\right)^{n-1}}{\left(\displaystyle\sum_{q=1}^{ Q}  \psi(q)\left(\log \psi(q)^{-1}\right)^{n-1}\right)^{-1}} \, >  \, 0  \ .
\end{equation}
In turn,
\begin{equation}\label{vb+12}
|\cS_n^\times(\psi)|   =  0 =  | \cD_n^\times(\psi)|
\qquad \text{if } \qquad  \sum_{q=1}^{\infty} \psi(q)\left(\log \psi(q)^{-1}\right)^{n-1}
  < \infty\ .
\end{equation}
\end{theorem}

\medskip

\noindent Note that the `additional' assumption  (\ref{cond1}) implies that the convergence/divergence properties of the sums within Conjectures \ref{con1} \& \ref{con2} are equivalent.

\medskip

\noindent{\em Remark 3.\ } Theorem \ref{MDST}  enables us to establish the complete  analogue of  Gallagher's  multiplicative  theorem \cite{gallagher1962} within  the framework of the `$p$-adic Littlewood Conjecture' -- see~\S\ref{king}. It is also worth pointing out that the same arguments that show that the Duffin-Schaeffer theorem is valid for example when  $\psi(q)$ is monotonic with $q$ restricted to a lacunary sequence, or  when $\psi(q) $ is arbitrary with $q$ restricted to the sequence of primes, are equally applicable within the context of Theorem \ref{MDST}.

\bigskip

\noindent{\em Remark 4.\ }
In the case when $\psi(q)\le q^{-\delta}$ for all sufficiently large $q\in\N$ and some fixed $\delta>0$ the term $\log \psi(q)^{-1}$ can replaced with $\log q$ throughout the statement of Theorem~\ref{MDST}.
However, in general this `modified' version (in which  $\log \psi(q)^{-1}$ is replaced with $\log q$) of the divergence part of Theorem~\ref{MDST}   is false\footnote{The authors are grateful to the anonymous reviewer of the paper who has pointed a mistake in the earlier version of the paper, where we miss out the fact that such a simplification is not always possible.}. For instance, assume that $n \geq 2$ and
let $P$ be an infinite collection of primes such that $\sum_{p\in P} (\log\log p)^{n-1}(\log p)^{-1}<\infty$ and
$\psi(q)=(\log q)^{-1}$ if $q \in P$ and $0$ otherwise. Then
$$
\sum_{q=1}^{ \infty}  \psi(q)\left(\log \psi(q)^{-1}\right)^{n-1} <\sum_{p\in P}^{ \infty}  \frac{(\log\log q)^{n-1}}{\log p}<\infty
$$
and so  $|\cS_n^\times(\psi)|=0=|\cD_n^\times(\psi)|$ by the convergence part of  Theorem~\ref{MDST}.  However,
$$
\sum_{q=1}^{ \infty}  \psi(q)(\log q)^{n-1}\ge \sum_{q=1}^{ \infty}  \frac{(\log q)^{n-1}}{\log q}\ge
\sum_{q=1}^{ \infty}  1=\infty
$$
since $ n \geq 2$, and
$$
\sum_{q=1}^{ Q}  \left(\frac{\varphi(q)}{q}\right)^n\psi(q)(\log q)^{n-1}\asymp
\sum_{q=1}^{ Q}  \psi(q)(\log q)^{n-1}\
$$
since $\psi$ is supported on primes and thus $\varphi(q)=q-1$ whenever $\psi(q)\neq0$. Thus, the `modified' version of the divergence part of Theorem~\ref{MDST} would be false for this $\psi$. Here and elsewhere $ a \asymp b $ means that $a$ and $b$ are comparable, that is $a \ll b$ and $a\gg b$, where $\ll$ and $\gg$ are the Vinogradov symbols that indicate an inequality with an unspecified positive multiplicative constant.

\section{Cross Fibering Principle \label{cross}}

Let $X$ and $Y$ be two non-empty sets.
Let $S\subset X\times Y$. Given $x\in X$, the set
$$
S_x:=\{y:(x,y)\in S\}\subset Y
$$
will be called a \emph{fiber of $S$ through $x$}. Similarly, given $y\in Y$, the set
$$
S^y:=\{x:(x,y)\in S\}\subset X
$$
will be called a \emph{fiber of $S$ through $y$}.
Given a measure $\mu$ over $X$, we will say that $A\subset X$
is \emph{$\mu$-trivial} if $A$ is either null or full with respect to  $\mu$; that is
$$\mu(A)=0  \qquad {\rm or }  \qquad
\mu(X\setminus A)=0  \, . $$
It is an immediate consequence of Fubini's theorem (see below) that
\begin{equation}\label{vb10}
\text{$S$ is $\mu\times\nu$-trivial}\quad
\Longrightarrow
\quad\text{$\mu$-almost every fiber $S_x$ is $\nu$-trivial},
\end{equation}
and likewise
\begin{equation}\label{vb11}
\text{$S$ is $\mu\times\nu$-trivial}\quad
\Longrightarrow
\quad\text{$\nu$-almost every fiber $S^y$ is $\mu$-trivial}.
\end{equation}
Neither of these implications can be reversed in their own right.
However, if the right hand side statements are combined together then we actually have a criterion which we will refer to as the {\em cross fibering principle}.

\begin{theorem}\label{t2}
Let $\mu$ be a $\sigma$-finite measure over $X$, $\nu$ be
a $\sigma$-finite measure over $Y$ and
$S\subset X\times Y$ be a $\mu\times\nu$-measurable set. Then
\begin{equation}\label{vb1}
\text{$S$ is $\mu\times\nu$-trivial}\quad
\Longleftrightarrow
\quad\begin{array}{ll}
\text{$\mu$-almost every fiber $S_x$ is $\nu$-trivial}\,  \\[1ex]
~ \hspace*{15ex} \& \\[1ex]
\text{$\nu$-almost every fiber $S^y$ is $\mu$-trivial}\,.
\end{array}
\end{equation}
\end{theorem}

The proof of this theorem will make use of
the following form of Fubini's theorem which can be found in
\cite[pg.233]{Billingsley-95:MR1324786} and \cite[\S2.6.2]{Federer-69:MR0257325}.

\begin{FubiniTheorem}
Let $\mu$ be a $\sigma$-finite measure over $X$ and $\nu$ be
a $\sigma$-finite measure over $Y$.
Then $\mu\times\nu$ is a regular measure over $X\times Y$ such that
\begin{enumerate}
\item[{\rm(i)}]
If $A$ is a $\mu$-measurable set and $B$ is a $\nu$-measurable set then $A\times B$
is a $\mu\times\nu$-measurable set and
$$
(\mu\times\nu)(A\times B)=\mu(A)\cdot\nu(B)  \, .
$$
\item[{\rm(ii)}]
If $S$ is a $\mu\times\nu$-measurable set, then
$$
S^y\quad\text{is $\mu$-measurable for $\nu$-almost all $y$},
$$
$$
S_x\quad\text{is $\nu$-measurable for $\mu$-almost all $x$},
$$
the functions
\begin{equation}\label{f1}
X\to\overline\R \,:\, x\mapsto\nu(S_x)\qquad\text{and}\qquad Y\to\overline\R \,:\, y\mapsto\mu(S^y)
\end{equation}
are integrable and
\begin{equation}\label{f2}
(\mu\times\nu)(S)=\int\mu(S^y)d\nu=\int\nu(S_x)d\mu.
\end{equation}
\end{enumerate}
\end{FubiniTheorem}

\subsection{Proof of Theorem~\ref{t2}}

The measures $\mu$ and $\nu$
are $\sigma$-finite.  Thus, without loss of generality we can assume that the measures are
finite and indeed that they are  probability measures; that is
$$\mu(X)=1 = \nu(Y)  \, . $$

\noindent {{\em Necessity }($\Longrightarrow$).}  Without loss of generality,  we can assume that
$(\mu\times\nu)(S)=0$ since otherwise we can replace $S$ by  its complement $X\setminus S$.
Therefore, both the integrals appearing in  (\ref{f2}) vanish.
Note that the integrals themselves are obtained by integrating the  non-negative functions (\ref{f1}). The upshot is that these functions vanish almost everywhere with respect to the appropriate measures which in turn  implies the
right hand side of (\ref{vb1}).

\medskip

\noindent {{\em Sufficiency }($\Longleftarrow$).}
Let $\tilde X$ be the set of $x\in X$ such that
$S_x$ is $\nu$-measurable and trivial. Similarly,  let $\tilde Y$ be the set of $y\in Y$ such that $S^y$ is $\mu$-measurable and trivial.
In view of part (ii) of Fubini's theorem and the right hand side  of (\ref{vb1})
we have that both $\tilde X$ and $\tilde Y$ are sets of full measure; that is
$\mu(X\setminus \tilde X)=0$ and $\nu(Y\setminus \tilde Y)=0$. In particular,
$\tilde X$ is  $\mu$-measurable  and  $\tilde Y$ is  $\nu$-measurable.
Now partition $\tilde X$ and $\tilde Y$ into two disjoint subsets as follows:
$$
\begin{array}{lcl}
  X_0:=\{x\in \tilde X:\nu(S_x)=0\}, && Y_0:=\{y\in \tilde Y:\mu(S^y)=0\},  \\[2ex]
  X_1:=\tilde X\setminus X_0=\{x\in \tilde X:\nu(S_x)=1\}, &\qquad&
  Y_1:=\tilde Y\setminus Y_0=\{y\in \tilde Y:\nu(S^y)=1\}.
\end{array}
$$
Let ${\cal X}_A$ denote the characteristic function of a set $A$.
By definition and part (ii) of Fubini's theorem,  the functions (\ref{f1})
almost everywhere coincide with the functions ${\cal X}_{X_1}$ and ${\cal X}_{Y_1}$. Since the functions (\ref{f1}) are integrable,
the functions ${\cal X}_{X_1}$ and ${\cal X}_{Y_1}$ are also integrable and so it follows that the sets $X_1$ and $Y_1$ are respectively $\mu$ and $\nu$-measurable.
This together with the fact that  $\tilde X$ and $\tilde Y$ are respectively  $\mu$ and $\nu$-measurable, implies that
$X_0=\tilde X\setminus X_1$ is $\mu$-measurable and $Y_0=\tilde Y\setminus Y_1$ is  $\nu$-measurable.

Observe that $\mu(X_0)+\mu(X_1)=\mu(\tilde X)=1$ and $\nu(Y_0)+\nu(Y_1)=\nu(\tilde Y)=1$.
Let us assume  that  the sets $X_i$ and $Y_i$ are non-trivial.   In other words,
\begin{equation} \label{contro}
  0<\mu(X_i)<1  \quad  {\rm and }  \quad  0<\nu(Y_i)<1   \qquad {\rm for} \quad  i=0,1  \, .
\end{equation}
By part (i) of Fubini's theorem, the set $M:=X_0\times Y_1$ is
$\mu\times\nu$-measurable. Now consider the set $S\cap M$ and observe that $M^y= X_0$ if $y\in Y_1$ and $M^y=\emptyset$ otherwise. Therefore,
on using the first equality of (\ref{f2}) we obtain that
\begin{equation}\label{vb2}
(\mu\times\nu)(S\cap M)
=\int\mu(S^y\cap M^y)d\nu=\int\mu(S^y\cap X_0) \ {\cal X}_{Y_1}(y)d\nu .
\end{equation}
By definition, for $y\in Y_1$ the set $S^y$ is full in $X$ and thus
is full in $X_0$. As a consequence, we have that $\mu(S^y\cap X_0)=\mu(X_0)$ for $y\in Y_1$. Therefore,
 \eqref{contro} and (\ref{vb2})  imply that
\begin{equation}\label{vb3}
(\mu\times\nu)(S\cap M)
=\int\mu(X_0){\cal X}_{Y_1}(y)d\nu=\mu(X_0)\nu(Y_1)>0.
\end{equation}
On the other hand, observe that  $M_x= Y_1$ if $x\in X_0$ and $M_x=\emptyset$ otherwise.
Then, on  using the second equality of (\ref{f2}) we obtain that
\begin{equation}\label{vb4}
(\mu\times\nu)(S\cap M)
=\int\nu(S_x\cap M_x)d\mu=\int\nu(S_x\cap Y_1){\cal X}_{X_0}(x)d\mu.
\end{equation}
By definition, for $x\in X_0$ the set $S_x$ is null and so
$\nu(S_x\cap Y_1)=0$ for $x\in X_0$. Therefore,
(\ref{vb4}) implies that
$$
(\mu\times\nu)(S\cap M)
=\int0 \, d\mu=0.
$$
This contradicts (\ref{vb3}). Therefore
at least one of the sets $X_i$ and $Y_i$ must be trivial.  This together with
(\ref{f2}) implies  that $S$ is trivial and thereby completes the proof.

\bigskip

\noindent{\em Remark 5.\ } Using induction Theorem~\ref{t2} can be easily extended to the product of any finite number of measure spaces.


%

\section{Proof of Theorem~\ref{t1}}

The proof is by induction. Consider the set $\cS_n^\times(\psi)$.
When $n=1$, we have that   $\cS_1^\times(\psi)=\cS_1(\psi)$ and Cassels' zero-one law   implies that  $\cS_1^\times(\psi)$ is $\mu$-trivial where $\mu$ is one-dimensional  Lebesgue measure on $X:=[0,1]$ .

Now assume that $n>1$ and that  Theorem~\ref{t1} is true for all dimensions $k <n$.
Given a $k$-tuple $(x_1,\dots,x_k)\in[0,1]^k$, consider  the function
$$
\psi_{(x_1,\dots,x_k)}(q) \, := \, \frac{\psi(q)}{\|qx_1\|\dots\|qx_k\|} \, .
$$
Here we adopt the convention  that $\alpha/0:=+\infty$ if $\alpha>0$ and that $\alpha/0:=0$ if $\alpha=0$. With reference to \S\ref{cross},
let $Y:=[0,1]^{n-1}$ and let $\nu$ be $(n-1)$-dimensional  Lebesgue measure on
 $Y$. Furthermore, let $S:=\cS^\times_n(\psi)$. Then it is
readily verified that for any $x_1\in X$ the fiber $S_{x_1}$ is equal to the set $\cS_1^\times(\psi_{(x_1)})$ and similarly for any
$(x_2,\dots,x_n)\in Y$ the fiber
$S^{(x_2,\dots,x_n)}$ is equal to the set  $\cS_{n-1}^\times(\psi_{(x_2,\dots,x_n)})$.
In view of the induction hypothesis, we have that   $S_{x_1}$ is $\mu$-trivial and $S^{(x_2,\dots,x_n)}$ is $\nu$-trivial.
Therefore, by Theorem \ref{t2} it follows that  $S$ is $ \mu \times \nu $-trivial. In other words,  the
$n$-dimensional Lebesgue measure of  $\cS^\times_n(\psi)$ is either zero or one. This establishes Theorem~\ref{t1}  for the set $\cS_n^\times(\psi)$.

\medskip

Apart from obvious notational changes, the proof for the set  $\cD^\times_n(\psi)$ is exactly the   same as above  except for that fact that when $n=1$  we appeal to  Gallagher's  zero-one law rather than Cassels' zero-one law.

\subsection{A multiplicative zero-one law for linear forms}

In what follows $m \geq1 $ and $n \geq 1 $ are integers.
 Given a `multi-variable' approximating function
$\Psi:\Z^n\to\Rp$, let $\cS^\times_{n,m}(\Psi)$ denote  the set of $\vv X\in[0,1]^{mn}$
such that
\begin{equation}\label{e:001}
 \Pi (\vv q\vv X+\vv p)<\Psi(\vv q)
\end{equation}

\noindent holds for infinitely many  $(\vv p,\vv q)\in\Z^n\times\Z^m\nz$.
Here $\Pi(\vv y):=\prod_{i=1}^n|y_i|$ for a vector  $\vv y=(y_1,\dots,y_n)\in\R^n$,
$\vv X$ is
regarded as an $m\times n$ matrix and $\vv q$ is regarded as a
row vector. Thus, $\vv q\vv X\in\R^n$ represents a system of $n$ real
linear forms in $m$ variables.   Naturally, let  $\cD^\times_{m,n}(\Psi)$ denote the subset of  $\cS^\times_{m,n}(\Psi)$ corresponding to $\vv X\in[0,1]^{mn}$ for which  \eqref{e:001} holds infinitely often with the additional co-primeness condition
$(p_i,\vv q)=1  $ for all $  1 \le i \le n$.
 Clearly, when $m=1$ and $\Psi(q) = \psi(|q|) $ the sets
$\cS^\times_{m,n}(\Psi)$ and $\cS^\times_{n}(\psi)$  coincide as do the sets $\cD^\times_{m,n}(\Psi)$ and $\cD^\times_{n}(\psi)$.

The following statement is the  natural generalisation of Theorem~\ref{t1} to the linear forms framework.  It also gives a positive answer to  Question~4 raised in \cite{Beresnevich-Velani-08:MR2457266}.

\begin{theorem}\label{t5} Let $m,n \geq1 $  and $\Psi:\Z^n\to\Rp $ be a non-negative  function. Then
$$|\cS^\times_{m,n}(\Psi)|\in\{0,1\}  \qquad  and   \qquad |\cD^\times_{m,n}(\Psi)|\in\{0,1\}  \, . $$
\end{theorem}

\medskip

\noindent In view of the linear forms version of the Cassels-Gallagher  zero-one law established
in \cite{Beresnevich-Velani-08:MR2457266}, the proof of  Theorem~\ref{t5} is pretty much the same as the proof of Theorem~\ref{t1} with obvious modification. More specifically,  all that is required from \cite{Beresnevich-Velani-08:MR2457266} is Theorem~1 with $n=1$.

\section{Proof of Theorem \ref{MDST} \label{ohya}}

To begin with, we recall that $ \cS_n^\times(\psi) \supset   \cD_n^\times(\psi) $   and therefore is suffices to prove the divergence part for $ \cD_n^\times(\psi)$ and the convergence part for $\cS^\times(\psi)$ only.
Regarding the divergence case, by Theorem~\ref{t1}, we are done if we can show that
$| \cD_n^\times(\psi)  |  > 0$. Given $q \in \N$, let
$$
A_q:=\big\{(x_1,\dots,x_n)\in[0,1)^n:x_1\cdots x_n\le \psi(q)\big\}
$$
and
\begin{eqnarray*}
B_q    :=  \left\{\vv x\in[0,1)^n\,:\begin{array}{l}
\, q\vv x-\vv p\in A_q\text{  for  some  $\vv p =(p_1, \ldots, p_n) \in \Z^n$}\\[1ex]
\text{with $(p_i,q)=1$ for all $i$}                                   \end{array}
\right\}.
\end{eqnarray*}
Note that $A_q =B_q =\emptyset $ whenever $\psi(q) =0$ and   that
$$B:=\limsup_{q \to \infty } B_q  \, \subset \, \cD_n^\times(\psi) \, . $$
Thus it suffices to prove that $|B|>0$.
For this purpose we will use the following generalisation of the divergent part of the
standard Borel-Cantelli lemma, see for example \cite[Lemma 5]{Sprindzuk-1979-Metrical-theory}.

\begin{lemma}\label{l4}
Let $(\Omega,A,\mu) $ be a probability  space and $\{E_q\} \subseteq A$ be a
sequence of sets such that $\sum_{q=1}^\infty \mu(E_q)=\infty $.
Then
$$ \mu( \limsup_{q \to \infty} E_q ) \; \geq \;  \limsup_{Q \to
\infty} \frac{ \left( \sum_{s=1}^{Q} \mu(E_s) \right)^2 }{ \sum_{s,
t = 1}^{Q} \mu(E_s \cap E_t ) }  \  \  \ . $$
\end{lemma}

\medskip

\noindent Naturally we shall use this lemma with $E_q=B_q$. The following estimates for the measure of $|B_q|$ can be found in \cite[\S\S1,2]{Gallagher-61:MR0133297}  -- they make use of the assumption that $0\le \psi(q)\le 1/2$.  For $ q \in \N$
$$
|A_q|\asymp \psi(q)\left(\log \psi(q)^{-1}\right)^{n-1}
$$
and
$$
|B_q|=(\varphi(q)/q)^n|A_q|\asymp (\varphi(q)/q)^n\psi(q)\left(\log \psi(q)^{-1}\right)^{n-1}\ .
$$
Then, by \eqref{cond1}, we have that for infinitely many $Q$
\begin{equation}\label{vb+1}
\sum_{q=1}^Q | B_q |  \asymp  \sum_{q=1}^Q \psi(q)\left(\log \psi(q)^{-1}\right)^{n-1}\asymp \sum_{q=1}^Q | A_q |  \, .
\end{equation}
Together with the  divergent sum hypothesis this implies that
\begin{equation} \label{obvious}
\sum_{q=1}^\infty | B_q |  = \infty \, .
\end{equation} Regarding the measures of overlaps, Lemma 2  in \cite{gallagher1962}  implies that
\begin{equation}\label{vb+2}
|B_q\cap B_{q'}|\le |A_q|\,|A_{q'}|\qquad\text{for }q\neq q'
\end{equation}
with $\psi(q)\neq0$ and $\psi(q')\neq0$. Note that \eqref{vb+2} is valid if $\psi(q)=0$ or $\psi(q')=0$ since we have zero on both sides of the inequality. Since \eqref{vb+1} diverges, $\sum_{q=1}^{ Q } |A_q|\ge1$ for $Q$ sufficiently large and so
$$
\sum_{q,q' =1}^{Q} | B_q\cap  B_{q'} |  \ \le \ \Big(\sum_{q=1}^{ Q } |A_q|\Big)^2 +\sum_{q=1}^{ Q } |A_q|\ \le \
2\Big(\sum_{q=1}^{ Q } |A_q|\Big)^2 .
$$
This together with \eqref{vb+1} and \eqref{obvious},  gives via Lemma~\ref{l4} that  $|\limsup_{q\to\infty}B_q|>0$ and thereby proves the divergence case of Theorem \ref{MDST}.

The convergence case is a consequence of, for example, Theorem~13 from \cite[\S5]{Sprindzuk-1979-Metrical-theory}. Before saying how to derive it note that the role of $m$ and $n$ is reversed in \cite[\S5]{Sprindzuk-1979-Metrical-theory} compared to the present paper. Thus, with reference to Theorem~13 in \cite[\S5]{Sprindzuk-1979-Metrical-theory} one has to take $n=1$, $S=\N$ and $A(\overline a)$ to consist of $(x_1,\dots,x_m)\in[0,1)^m$ such that $\|x_1\|\cdots\|x_m\|<\psi(\overline a)$ for $\overline a\in S$. Note that any point in $A(\overline a)$ is obtained from a point of $A_q$ (introduced above) by applying relevant symmetries $x_i\mapsto 1-x_i$. This gives that $|A(\overline a)|\ll |A_q|\ll \psi(q)\left(\log \psi(q)^{-1}\right)^{n-1}$ and hence the condition $\sum_{\overline a}|A(\overline a)|<\infty$ which is required to derive our Theorem~\ref{MDST} for convergence from Theorem~13 in \cite[\S5]{Sprindzuk-1979-Metrical-theory}.

\subsection{An application to $p$-adic approximation \label{king}}

Theorems \ref{t1} \& \ref{MDST}  settle the conjecture and problem stated in \cite[\S4.5]{bhv} regarding  the multiplicative set $ \cS_n^\times(\psi)$.     In particular, as a consequence of Theorem \ref{MDST} we are able to prove the following generalisation of the main result appearing in \cite{bhv}.  In short the statement   corresponds to the complete  analogue of  Gallagher's  multiplicative  theorem \cite{gallagher1962} within the framework of the `$p$-adic Littlewood Conjecture' -- for further details see \cite{badvel,bhv} and references within.  Given a prime $p$, we denote by $|q|_p $ the $p$-adic norm of $q \in \Z$.

\begin{theorem}
Let $k\in\N$, $p_1,\ldots, p_k$ be distinct prime numbers and $F:\N\to\Rp$ be a positive function such that
\begin{equation}\label{F}
\text{$F(q)=F(q')$ whenever $|q|_{p_i}=|q'|_{p_i}$ for all $i$.}
\end{equation}
Let $\Psi:\N\to\Rp$ will be a positive decreasing function. If
\begin{equation}\label{vb+8}
\sum_{q=1}^{\infty}  \;  \frac{\Psi(q)}{F(q)}  \ \left(\log_+\frac{F(q)}{\Psi(q)}\right)^{n-1}
\end{equation}
converges, where $\log_+x:=\log\max\{2,x\}$, then for almost every $(x_1,\ldots ,x_n)\in \R^n $ the inequality
\begin{equation}\label{vb+8++}
F(q)\|qx_1\|\cdots\|qx_n\|<\Psi (q)
\end{equation}
has only finitely many solutions $q \in \N$. On the other hand, if \eqref{vb+8} diverges then for almost every $(x_1,\ldots ,x_n)\in \R^n $ inequality \eqref{vb+8++} has infinitely many solutions $q \in \N$.
\end{theorem}

\bigskip

\noindent{\em Proof of Theorem~5.} We will use Theorem~\ref{MDST} with $\psi(q)=\Psi(q)/F(q)$. Indeed, the case of convergence is a straightforward application of Theorem~\ref{MDST} as in this case the convergence of \eqref{vb+8} implies the convergence condition within \eqref{vb+12}.

In what follows we consider the divergence case. First of all observe that $\|qx_1\|\cdots\|qx_n\|\le 2^{-n}$ for all $q\in\N$. Therefore, without loss of generality we can assume that $\Psi(q)/F(q)<2^{-n}$ for all $q\in\N$. Furthermore, by replacing $\Psi(q)$ with $2^{-n}\Psi(q)$ if necessary, we can assume that $\Psi(q)/F(q)<e^{-n}$ for all $q$. In particular, this means that $\log_+\frac{F(q)}{\Psi(q)}=\log\frac{F(q)}{\Psi(q)}$.

\noindent Throughout, $\Z_+$ will denote non-negative integers, $\alpha=(\alpha_1,\dots,\alpha_k)\in\Z_+^k$
and $|\alpha|=\max_i\alpha_i$. Each $q\in\N$ can be uniquely written as $q=p_1^{\alpha_1}\dots p_k^{\alpha_k}\ell$ for some $\alpha=(\alpha_1,\dots,\alpha_k)$ and $\ell\in\N$ with $(\ell,p_1\cdots p_k)=1$. Then, by \eqref{F}, the monotonicity of $\Psi$ and the assumption that  $\Psi(q)/F(q)<e^{-n}$, the function
$$
f_\alpha(\ell):=\frac{\Psi(q)}{F(q)}  \ \left(\log\frac{F(q)}{\Psi(q)}\right)^{n-1}\quad\text{where }q=p_1^{\alpha_1}\dots p_k^{\alpha_k}\ell\quad\text{and}\quad (\ell,p_1\cdots p_k)=1
$$
is monotonically decreasing in $\ell$ for each fixed $\alpha$. Label the numbers $\ell_i$ ($i\in\N$) with $(\ell_i,p_1\cdots p_k)=1$ in increasing order; i.e.  $\ell_1 < \ell_2< \ell_3 < \dots $. By \eqref{vb+8},
\begin{equation}\label{vb+4}
\sum_{\alpha\in\Z_+^k}\ \sum_{i=1}^\infty \;  f_\alpha(\ell_i)
=\infty.
\end{equation}
Thus, by Theorem~\ref{MDST} with $\psi(q)=\Psi(q)/F(q)$, to complete the proof of Theorem~5 it suffices to show that for sufficiently large $Q$
\begin{equation}\label{vb+9}
    \sum_{|\alpha|\le Q}\ \sum_{i\le Q} \;  \left(\frac{\varphi(p_1^{\alpha_1}\dots p_k^{\alpha_k}\ell_i)}{p_1^{\alpha_1}\dots p_k^{\alpha_k}\ell_i}\right)^nf_\alpha(\ell_i)
    \gg     \sum_{|\alpha|\le Q}\ \sum_{i\le Q} \;  f_\alpha(\ell_i)  \, .
\end{equation}
Since $\varphi(p_1^{\alpha_1}\cdots p_k^{\alpha_k}\ell_i)= \prod_{i=1}^k(1-p_i^{-1})p_1^{\alpha_1}\cdots p_k^{\alpha_k}\varphi(\ell_i)$, inequality \eqref{vb+9} would follow  on showing that for each fixed $\alpha\in\Z_+^k$
\begin{equation}\label{vb+10}
 \sum_{i\le Q}  \;  \left(\frac{\varphi(\ell_i)}{\ell_i}\right)^nf_\alpha(\ell_i)
    \gg      \sum_{i\le Q}  \;  f_\alpha(\ell_i)
\end{equation}
with the implied constant being independent of $\alpha$.
Lemma~2 from \cite{bhv} gives  that  $\sum_{j\le i}\varphi(\ell_j)/\ell_j\gg i$. This together with Jensen's inequality implies  that
$\sum_{j\le i}\big(\varphi(\ell_j)/\ell_j\big)^n\gg i$.
Then, by  partial summation  and the monotonicity of $f_\alpha$, for each fixed $\alpha$ and $Q>1$ we have that
$$
\sum_{i\le Q}\left(\frac{\varphi(\ell_i)}{\ell_i}\right)^n f_\alpha(\ell_i)=\sum_{i\le Q}(f_\alpha(\ell_i)-f_\alpha(\ell_{i+1}))\sum_{j=1}^i\left(\frac{\varphi(\ell_j)}{\ell_j}\right)^n +
$$

$$
\qquad \qquad \qquad \qquad \qquad \qquad \qquad \qquad +f_{\alpha}(\ell_{Q+1})\sum_{j=1}^Q\left(\frac{\varphi(\ell_j)}{\ell_j}\right)^n
$$

$$
\qquad \qquad \gg \sum_{i\le Q}i\big(f_\alpha(\ell_i)-f_\alpha(\ell_{i+1})\big) +Qf_{\alpha}(\ell_{Q+1})=
\sum_{i\le Q}f_\alpha(\ell_i).
$$
This establishes \eqref{vb+10} and thus completes the proof.

\bigskip
\bigskip

\noindent{\em Remark 6.\ } It is impossible to replace $\log\frac{F(q)}{\Psi(q)}$ with $\log q$ within \eqref{vb+8}.
To see that this is so, let $k=n=1$, $p_1=p$, $\Psi(q)=q^{-2}$ and $F(q)=|q|_p^2(1-\log|q|_p)^2$.
Write each $q\in\N$ as $p^\alpha\ell$ with $\alpha\in\Z_+$, $\ell\in\N$ and $(p,\ell)=1$. Then,
$$
\frac{\Psi(q)}{F(q)}=\frac{(p^\alpha\ell)^{-2}}{p^{-2\alpha}(1+\alpha\log p)^{2}} \asymp\frac{1}{\ell^2(1+\alpha)^{2}}.
$$
Consequently \eqref{vb+8} is comparable to
$$
\sum_{\alpha\ge0}\ \sum_{(\ell,p)=1}  \;  \frac{1}{\ell^2(1+\alpha)^2}  \ \log (\ell^2(1+\alpha)^2) \ll
\sum_{\alpha\ge0}\ \sum_{\ell\ge 1}  \;  \frac{(\log\ell)(\log(1+\alpha))}{\ell^2(1+\alpha)^2}  <\infty\,.
$$
On the other hand,
$$
\sum_{q=1}^{\infty}  \;  \frac{\Psi(q)}{F(q)}  \ \left(\log q\right)^{n-1} \asymp \sum_{\alpha\ge 0}\ \sum_{(\ell,p)=1}  \;  \frac{1}{\ell^{2}(1+\alpha)^2}  \ \log (p^\alpha\ell) \gg
\sum_{\alpha\ge 0}\;  \frac{1}{1+\alpha}  =\infty\,.
$$

\vspace*{5ex}

\vspace*{2ex}

{\small
\noindent Victor V. Beresnevich: Department of Mathematics,
University of York,

\vspace{0mm}

\noindent\phantom{Victor V. Beresnevich: }Heslington, York, YO10
5DD, England.

\vspace{0mm}
\noindent\phantom{Victor V. Beresnevich: }e-mail: vb8@york.ac.uk

\vspace{3mm}

\noindent Alan K. Haynes: Department of Mathematics, University of
York,

\vspace{0mm}

\noindent\phantom{Alan K. Haynes: }Heslington, York, YO10 5DD,
England.


\noindent\phantom{Alan K. Haynes: }e-mail: akh502@york.ac.uk

\vspace{3mm}

\noindent Sanju L. Velani: Department of Mathematics, University of
York,

\vspace{0mm}

\noindent\phantom{Sanju L. Velani: }Heslington, York, YO10 5DD,
England.


\noindent\phantom{Sanju L. Velani: }e-mail: slv3@york.ac.uk

}

\end{document}